\documentclass[11pt,reqno]{amsart}

\theoremstyle{plain}
\newtheorem{theorem} {Theorem}[section]
\newtheorem{lemma}[theorem] {Lemma}

\theoremstyle{definition}

\theoremstyle{remark}
\newtheorem{remark}[theorem] {Remark}

\numberwithin{equation}{section}

\usepackage{amsfonts}
\usepackage{amssymb}
\usepackage{amscd}

\newcommand{\R}{{\mathbb R}}

\newcommand{\N}{{\mathbb N}}

\newcommand{\PP}{{\mathcal P}}
\newcommand{\F}{{\mathcal F}}

\newcommand{\CC}{{\mathbb C}}

\newcommand{\E}{{\mathcal E}}
\newcommand{\al}{{\alpha}}
\newcommand{\la}{{\lambda}}
\newcommand{\sa}{{\sigma}}

\newcommand{\iy}{{\infty}}
\newcommand{\vphi}{{\varphi}}
\newcommand{\vep}{{\varepsilon}}
\newcommand{\g}{{\gamma}}
\newcommand{\de}{{\delta}}

\newcommand{\be}{{\beta}}

\newcommand{\bna}{\begin{eqnarray}}
\newcommand{\ena}{\end{eqnarray}}
\newcommand{\ba}{\begin{eqnarray*}}
\newcommand{\ea}{\end{eqnarray*}}
\newcommand{\beq}{\begin{equation}}
\newcommand{\eeq}{\end{equation}}

\begin{document}

\title[Asymptotic Behaviour of the Error]
{Asymptotic Behaviour of the Error of Polynomial Approximation
of Functions Like $\vert x\vert^{\al+i\be}$}
\author{Michael I. Ganzburg}
 \address{212 Woodburn Drive\\ Hampton,
 VA 23664\\USA}
 \email{michael.ganzburg@gmail.com}
 \begin{abstract}
 New asymptotic relations between the
$L_p$-errors of best approximation of univariate functions
by algebraic polynomials and entire functions of exponential type
are obtained for $p\in (0,\iy]$.
General asymptotic relations are applied to functions
$\vert x\vert^{\al+i\be},\,\vert x\vert^{\al}\cos(\be\log\vert x\vert)$,
and \linebreak $\vert x\vert^{\al}\sin(\be\log\vert x\vert)$.

 \end{abstract}
 \maketitle

 \section{Introduction}\label{S1}
 \noindent
\setcounter{equation}{0}
In this paper we discuss asymptotic relations between the
errors of best approximation of univariate functions
in the $L_p$-metrics, $p\in (0,\iy]$,
by algebraic polynomials and entire functions of exponential type.
General asymptotic relations are applied to functions
$\vert x\vert^{\al+i\be},\,\vert x\vert^{\al}\cos(\be\log\vert x\vert)$,
and $\vert x\vert^{\al}\sin(\be\log\vert x\vert)$.\vspace{.1in}\\
\textbf{Notation and Preliminaries.}
Let $\R$ be the set of all real numbers,
$\CC$ be the set of all complex numbers,
$\PP_n$ be the set
 of all algebraic polynomials with complex coefficients
 of degree at most $n,\,n\in \N:=\{1,\,2,\ldots,\}$,
 and $B_\sa$ be
 the set of all complex-valued entire
 functions of exponential type $\sa \ge 0$.
 Let $M_{K,N}$
 be the space of all measurable functions $f:\R\to\CC$
 of power growth, that is,
 \ba
 \vert f(x)\vert\le K(1+\vert x\vert)^N,\qquad x\in\R,
 \ea
 where $N\ge 0$ and $K>0$ are constants independent of $x$.
 Let $L_p(\Omega)$ be the space of all measurable complex-valued functions $f$
 on a measurable
 subset $\Omega$ of $\R$ with the finite quasinorm
 \beq\label{E1.1}
 \|f\|_{L_p(\Omega)}:=\left\{\begin{array}{ll}
 \left(\int_\Omega\vert f(x)\vert^p dx\right)^{1/p}, & 0<p<\iy,\\
 \mbox{ess} \sup_{x\in \Omega} \vert f(x)\vert, &p=\iy.
 \end{array}\right.
 \eeq
 Next, we define two errors of best approximation of $f$
 in the $L_p$-metrics,
 $p\in (0,\iy]$, by
\beq\label{E1.2}
 E_n\left(f,L_p([-a,a])\right):=\inf_{P\in\PP_n}\| f-P\|_{L_p([-a,a])},
 \quad f\in L_p([-a,a]);
\qquad A_\sa(f)_p:=\inf_{g\in B_\sa}\| f-g\|_{L_p(\R)},
\eeq
where $n\in\N,\,a>0$, and $\sa>0$.
Note that in the definition of $ A_\sa(f)_p,\,f$ is a measurable function on $\R$
but $f$ does not necessarily belong to $L_p(\R)$. In particular, $f$ could belong to
$M_{K,N}$
or be just locally integrable on $\R$. Let us also set
$A_{\sa\pm 0}(f)_p=\lim_{\tau\to \sa\pm}A_{\tau}(f)_p$.
In addition, note that the quasinorm $\|\cdot\|_{L_p(\Omega)}$
 in \eqref{E1.1}
and the quasiseminorm $E_n(\cdot,L_p([-a,a]))$ in \eqref{E1.2}
allow the following "triangle" inequalities:
 \bna
 \left\|f+g\right\|^{\tilde{p}}_{L_p(\Omega)}
 &\le& \left\| f\right\|^{\tilde{p}}_{L_p(\Omega)}
 + \left\| g\right\|^{\tilde{p}}_{L_p(\Omega)},
 \qquad f,\,g\in L_p(\Omega); \label{E1.3}\\
 E_n^{\tilde{p}}\left(f+g,L_p([-a,a])\right)
 &\le& E_n^{\tilde{p}}\left(f,L_p([-a,a])\right)
 +E_n^{\tilde{p}}\left(g,L_p([-a,a])\right),
 \quad f,\,g\in L_p([-a,a]),\label{E1.4}
 \ena
 where $\tilde{p}:=\min\{1,p\}$ for $p\in(0,\iy]$.
 We also need the following simple property
 ($f\in L_p([-a,a]), \linebreak p\in(0,\iy]$):
 \beq\label{E1.4a}
 E_n\left(f,L_p([-a,a])\right)
 = \vert\eta\vert^{1/p}E_n\left(f(\eta\cdot),
 L_p([-a/\vert\eta\vert,a/\vert\eta\vert])\right),
 \qquad \eta\in\R\setminus\{0\}.
 \eeq
 In addition, we use the notation
 \beq\label{E1.5}
 f_{\al,\be}(x):=\vert x\vert^{\al+i\be}
 =f_{\al,\be,c}(x)+if_{\al,\be,s}(x)
 :=\vert x\vert^{\al}\cos(\be\log\vert x\vert)
+i\vert x\vert^{\al}\sin(\be\log\vert x\vert).
\eeq

The Fourier transform of a function $\vphi\in L_1(\R)$
is defined by the formula
 \ba
 \F(\vphi)(x):=(2\pi)^{-1/2}\int_\R \vphi(t)\exp[-ixt]dt.
 \ea
 We use the same notation $\F(\varPhi)$ for the Fourier transform of a
 tempered distribution
 $\varPhi$ on $\R$.
 By definition (see, e.g., \cite[Sect. 4.1]{S1994}),
 $\varPhi$ is a continuous linear functional $\langle\varPhi,\vphi\rangle$
 on the Schwartz
 class $S(\R)$ of all
 test functions $\vphi$ on $\R$, and $\F(\varPhi)$ is defined by the formula
 $\langle \F(\varPhi),\vphi\rangle:=\langle\varPhi,\F(\vphi)\rangle,\,\vphi\in S(\R).$
 The Fourier transform inversion formula  is $\F(\F(\varPhi))=\Check{\varPhi}$,
 where $\langle \Check{\varPhi},\vphi\rangle:=\langle\varPhi,\vphi(-\cdot)\rangle,\,
 \vphi\in S(\R)$.
 In some cases
 discussed in the paper,
 functionals $\langle\varPhi,\cdot\rangle$ or $\langle \F(\varPhi),\cdot\rangle$ are
 generated by integrable or locally integrable functions.
 In particular, if $\varPhi$ is a locally
 integrable function on $\R$ and $\varPhi\in M_{K,N}$, then
$\langle \varPhi,\vphi\rangle
:=\int_\R \varPhi(t)\vphi(t)dt$ for all $\vphi\in S(\R)$.
If $\F(\varPhi)$ is a tempered distribution on $\R$ whose restriction $h$
 to $\R\setminus (-a,a),\,a>0$, is an integrable function, then
 $\langle\F(\varPhi),\vphi\rangle:=\int_\R h(t)\vphi(t)dt$ for all $\vphi\in S(\R)$ with
 $\mbox{supp}\, \vphi\subseteq \R\setminus (-a,a)$.

 Throughout the paper $C,\,C_1,\,C_2,\ldots$ denote positive
constants independent
of essential parameters.
 Occasionally we indicate dependence on certain parameters.
 The same symbol $C$ does not
 necessarily denote the same constant in different occurrences.
In addition, we use the floor and ceiling functions
 $\lfloor a \rfloor$ and $\lceil a \rceil$
and also use the generic notation $\Gamma(z)$ for the gamma function. \vspace{.12in}\\
\noindent
\textbf{Limit Relations between $E_n(f,L_p([-a,a]))$ and $A_\sa(f)_p$.}
It was Bernstein who in 1913 and 1938 initiated the study of asymptotic relations in
approximation theory
 by proving the following celebrated result:
 for any $\al>0$, there is a constant
$\mu_\al\in(0,\iy)$ such that
\beq\label{E1.6}
\lim_{n\to\iy}n^\al E_n\left(f_{\al,0},L_\iy([-1,1])\right)=\mu_\al.
\eeq
The proof of \eqref{E1.6} for $\al=1$ in \cite{B1913} was very long
 and very difficult.
 The case of $\al>0$ was discussed
 25 years later
 in a much shorter publication \cite{B1938},
 and the proof of \eqref{E1.6} was based on the limit equality
\beq\label{E1.7}
\lim_{n\to\iy}n^\al E_n\left(f_{\al,0},L_\iy([-1,1])\right)
=A_1\left(f_{\al,0}\right)_\iy.
\eeq
However, the proof of \eqref{E1.7} was still difficult.
In 1946 Bernstein \cite{B1946} came up with an impressive
general version of
\eqref{E1.7} for any function $f\in M_{K,N}$  and $\sa>0$:
\beq\label{E1.8}
A_\sa(f)_\iy\le \liminf_{n\to\iy}E_n\left(f,L_\iy([-n/\sa,n/\sa])\right)
\le \limsup_{n\to\iy}E_n\left(f,L_\iy([-n/\sa,n/\sa])\right)
\le A_{\sa-0}(f)_\iy.
\eeq
Later Bernstein \cite{B1947} proved that \eqref{E1.8} is valid for any
measurable function $f$ with $A_{\sa_0}(f)_\iy<\iy$ and $\sa>\mu \sa_0>0$,
where $\mu=1.508879\ldots$ is the positive solution to the equation
\beq\label{E1.8a}
\sqrt{x^2+1}/x=\log\left( \sqrt{x^2+1}+x\right).
\eeq
The detailed proof  of \eqref{E1.8}
with both versions of conditions on $f$ can be found in
\cite[Sects. 2.6.22 and 5.4.6]{T1963} and
  \cite[Appendix, Sect. 84]{A1965}.
It immediately follows from \eqref{E1.8} that for a. a. $\sa>0$
and $f\in M_{K,N}$,
\beq\label{E1.9}
 \lim_{n\to\iy}E_n\left(f,L_\iy([-n/\sa,n/\sa])\right)
= A_{\sa}(f)_\iy.
\eeq
For certain functions from $M_{K,N}$ equality \eqref{E1.9} holds for all $
\sa>0$. In particular, \eqref{E1.9} is valid for $f=f_{\al,0},\,\al>0$,
and all $\sa>0$. This fact
(along with \eqref{E1.4a}) immediately implies \eqref{E1.7}
(see also \cite{B1946.1}),
and such a proof of \eqref{E1.7} is shorter and more straightforward
than that in \cite{B1938}.

The author  \cite{G1991} showed that \eqref{E1.9} does not hold for
$f(x)=\cos x$ and $\sa=1$, but the following version of \eqref{E1.9}
is still valid for $f\in M_{K,N}$ and all $\sa>0$:
\beq\label{E1.10}
 \lim_{n\to\iy}E_n\left(f,L_\iy([-(n-\g_n)/\sa,(n-\g_n)/\sa])\right)
= A_{\sa}(f)_\iy,
\eeq
where $\g_n\in [0,n],\,n\in\N,\,
\liminf_{n\to\iy}(\g_n n^{-1/3}\log^{-2/3} n)>1/2$,
and $\lim_{n\to\iy}\g_n/n=0$.
In a sense, these conditions cannot be improved.

Raitsin \cite[Theorem B]{R1968} extended \eqref{E1.9} to all $\sa>0$
and to the $L_p$-metric in the form
\beq\label{E1.11}
 \lim_{n\to\iy}E_n\left(f,L_p([-n/\sa,n/\sa])\right)
= A_{\sa}(f)_p,\qquad p\in [1,\iy),
\eeq
provided that $f\in L_p(\R)$. An extension of \eqref{E1.11}
to functions with $f\in L_p(\R)$ or $f\in M_{K,N},\linebreak
A_{\sa_0}(f)_p<\iy,\,p\in(0,\iy),\,0<\sa_0<\sa$,
was obtained by the author \cite[Theorem 11.3]{G2000}.
Various extensions of relations \eqref{E1.8}, \eqref{E1.9},
\eqref{E1.10}, and \eqref{E1.11} to multivariate functions
 were given by the author \cite{G1982, G1991, G2000}.

 A detailed analysis of the relations
\beq\label{E1.12}
\lim_{n\to\iy}n^{\al+1/p} E_n\left(f_{\al,0},L_p([-1,1])\right)
=A_1\left(f_{\al,0}\right)_p
=\left\|f_{\al,0}-g_{\al,p}\right\|_{L_p(\R)},
\, p\in[1,\iy],\,\al>-1/p,
\eeq
was provided by Lubinsky \cite[Theorems 1.1,1.2,1.3]{L2007}.
In particular, various representations of the entire function
$g_{\al,p}\in B_1$ of best $L_p$-approximation to $f_{\al,0}$,
 and the uniqueness of $g_{\al,1}$ and
$g_{\al,\iy}$ were established in  \cite{L2007}.
The uniqueness of $g_{\al,p},\, p\in(1,\iy)$,
follows from a general result in strictly normed spaces
\cite [Sect. 9]{A1965}.
Limit relations like \eqref{E1.8}, \eqref{E1.11}, and
 \eqref{E1.12} in univariate spaces with exponential
 weights were discussed in \cite{G2008book}.

The Bernstein constants $A_1\left(f_{\al,0}\right)_p$
are known only for $p=1$ and $p=2$:
\bna
&&A_1\left(f_{\al,0}\right)_1=\frac{8\left\vert \sin \frac{\alpha \pi }{2}
\right\vert \Gamma \left( \alpha +1\right) }{\pi } \sum_{k=0}^{\iy}
\frac{ (-1) ^{k}}{( 2k+1)^{\alpha +2}},
\qquad \al>-1, \label{E1.13}\\
&&A_1\left(f_{\al,0}\right)_2=\frac{2\left\vert \sin \frac{\alpha \pi }{2}
\right\vert \Gamma \left( \alpha +1\right)}{ \sqrt{\pi \left( 2\alpha
+1\right) }},\qquad \al>-1/2.\label{E1.14}
\ena
A direct proof of \eqref{E1.13} was found by the author
\cite[Corollary 2]{G1983} (see also \cite[Example 3.2]{G2010});
certain special cases were discussed in \cite{L2005, CV2010}).
By a combined effort of Nikolskii \cite{N1947} and Bernstein
\cite{B1954}, it was shown that
$\lim_{n\to\iy}n^{\al+1} E_n\left(f_{\al,0},L_1([-1,1])\right)$
coincides with
the right-hand side of equality \eqref{E1.13}.

Equality \eqref{E1.14} was established by Raitsin \cite{R1969}.
The following extension of \eqref{E1.14} to $f_{\al,\be}$ was given
by the author and Lubinsky \cite[Eq. (2.8)]{GL2008}:
\beq\label{E1.15}
A_1\left(f_{\al,\be}\right)_2=\frac{2\left\vert \sin \frac{(\alpha+i\be) \pi }{2}
 \Gamma \left( \alpha+i\be +1\right)\right\vert}{ \sqrt{\pi \left( 2\alpha
+1\right) }},\qquad \al>-1/2,\quad \be\in\R.
\eeq
The estimates
\beq\label{E1.16}
C_1(\al,\be)n^{-\al}\le E_n(f,L_\iy([-1,1]))\le C_2(\al,\be)n^{-\al},
\qquad \al\in(0,1),\quad \be\in\R,\quad n\in\N,
\eeq
where  $f$ is one of the functions
$f_{\al,\be},\,f_{\al,\be,c}$, or $f_{\al,\be,s}$,
were proved by the author \cite[Lemma 2]{G2008}.
\vspace{.12in}\\
  \textbf{Main Results.}
  In this paper we prove more general versions of limit relations
\eqref{E1.8} for $p=\iy$ and
\eqref{E1.11} for $p\in(0,\iy)$ (Theorems \ref{T1.1}, \ref{T1.2},
and \ref{T1.2a}).
Next, we apply these results to functions
$f_{\al,\be},\,f_{\al,\be,c}$, and $f_{\al,\be,s}$,
defined by \eqref{E1.5},
generalizing asymptotic relations \eqref{E1.7} and \eqref{E1.12}
and improving estimates \eqref{E1.16} (Theorem \ref{T1.3}).
Finally, in addition to \eqref{E1.13} and \eqref{E1.15},
we find one more Bernstein constant $A_1(f_{0,\be})_\iy$ and also
associate constants (Theorem \ref{T1.4}).

 \begin{theorem} \label{T1.1}
 Let a sequence $\left\{a_n\right\}_{n=1}^\iy$
 of positive numbers
 satisfy the condition
 \beq\label{E1.16a}
 \sup_{n\in\N}\vert a_n- n\vert\le C<\iy.
 \eeq
 If $f\in M_{K,N}$ and $A_{\sa_0}(f)_p<\iy$ for $\sa_0>0$
 and $p\in(0,\iy]$, then for $\sa>\sa_0$ the following
 statements hold true:\\
 (a) If $p=\iy$, then
 \beq\label{E1.17}
A_\sa(f)_\iy
\le \liminf_{n\to\iy}E_n\left(f,L_\iy([-a_n/\sa,a_n/\sa])\right)
\le \limsup_{n\to\iy}E_n\left(f,L_\iy([-a_n/\sa,a_n/\sa])\right)
\le A_{\sa-0}(f)_\iy.
\eeq
(b) If $p\in(0,\iy)$, then
\beq\label{E1.18}
 \lim_{n\to\iy}E_n\left(f,L_p([-a_n/\sa,a_n/\sa])\right)
= A_{\sa}(f)_p.
\eeq
 \end{theorem}
 \noindent
 In the next two theorems the cases of $f\in L_p(\R)$
 and $A_{\sa_0}(f)_p<\iy$ are discussed.

\begin{theorem} \label{T1.2}
 Let  a sequence $\left\{a_n\right\}_{n=1}^\iy$
 of positive numbers
 satisfy condition \eqref{E1.16a} of Theorem \ref{T1.1}.
If $f\in L_p(\R),\,p\in(0,\iy]$, then for all $\sa>0$
statements (a) and (b) of Theorem \ref{T1.1} are valid.
\end{theorem}
\noindent
We recall that
$\mu=1.508879\ldots$ is the positive solution to equation
\eqref{E1.8a}.

\begin{theorem} \label{T1.2a}
 Let a sequence $\left\{a_n\right\}_{n=1}^\iy$
 of positive numbers
 satisfy condition \eqref{E1.16a} of Theorem \ref{T1.1}.
If $A_{\sa_0}(f)_p<\iy$ for $\sa_0>0$
 and $p\in(0,\iy]$, then for $\sa>\mu\sa_0$
statements (a) and (b) of Theorem \ref{T1.1} are valid.
\end{theorem}

Next, we apply Theorems \ref{T1.1} and \ref{T1.2a} to functions
$f_{\al,\be},\,f_{\al,\be,c}$, and $f_{\al,\be,s}$.

\begin{theorem} \label{T1.3}
If $\be\in\R$ and either $\al>\max\{-1,-1/p\},\,
p\in(0,\iy]$ or $\al=0,\,p=\iy$,
then the following
 statements hold true:\\
(a) If $f$ is one of the functions
$f_{\al,\be},\,f_{\al,\be,c}$, or $f_{\al,\be,s}$,
then $A_\sa(f)_p<\iy$ for $\sa>0$.\\
(b) The following equality holds true:
\beq\label{E1.19}
\lim_{n\to\iy}n^{\al+1/p} E_n\left(f_{\al,\be},L_p([-1,1])\right)
=A_1\left(f_{\al,\be}\right)_p.
\eeq
(c) For $\be\neq 0$ there are two sequences $\left\{n_k\right\}_{k=1}^\iy$
and $\left\{m_k\right\}_{k=1}^\iy$
of positive integers such that
\bna
\lim_{k\to\iy}n^{\al+1/p}_k E_{n_k}\left(f_{\al,\be,c},L_p([-1,1])\right)
&=& \lim_{k\to\iy}m^{\al+1/p}_k E_{m_k}\left(f_{\al,\be,s},L_p([-1,1])\right)
\nonumber\\
&=&A_1\left(f_{\al,\be,c}\right)_p,\label{E1.20}\\
\lim_{k\to\iy}n^{\al+1/p}_k E_{n_k}(f_{\al,\be,s},L_p([-1,1]))
&=& \lim_{k\to\iy}m^{\al+1/p}_k E_{m_k}\left(f_{\al,\be,c},L_p([-1,1])\right)
\nonumber\\
&=&A_1\left(f_{\al,\be,s}\right)_p.\label{E1.21}
\ena
\end{theorem}

New Bernstein constants are found in the following theorem.

 \begin{theorem} \label{T1.4}
 If $\be\in\R\setminus\{0\}$ and $\sa>0$, then
\beq\label{E1.22}
A_\sa\left(f_{0,\be}\right)_\iy
=A_\sa\left(f_{0,\be,c}\right)_\iy
=A_\sa\left(f_{0,\be,s}\right)_\iy=1,
\eeq
and entire functions of exponential type $\sa$
of best approximation to
$f_{\al,\be},\,f_{\al,\be,c}$, and $f_{\al,\be,s}$
in the metric of $L_\iy(\R)$ coincide with zero.
 \end{theorem}

 \begin{remark}\label{R1.5}
 Special cases of Theorems \ref{T1.1},
 \ref{T1.2}, and \ref{T1.2a} for
 $a_n=n$ (see \eqref{E1.8} and \eqref{E1.11})
    were established in \cite{B1946, B1947, R1968, G2000}.
    Note that relations \eqref{E1.17} and \eqref{E1.18}
     are used in the proof
    of Theorem \ref{T1.3} (c),
    and the replacement of $n$
     with a more general sequence
    $a_n,\,n\in\N$, in
    \eqref{E1.8} and \eqref{E1.11} is essential in this proof.
    A special case of Theorem \ref{T1.3} (b) for $\be=0$
    (see \eqref{E1.12}) was proved in \cite{L2007}.
\end{remark}

\begin{remark}\label{R1.6}
The corresponding versions of Theorems \ref{T1.3} and
\ref{T1.4} are also valid for the function
$\vert x\vert^{\al+i\be}(\mbox{sgn}\, x)$
 or the more general function
\ba
F_{\al,\be,a,b}(x):=\left\{\begin{array}{ll}
ax^{\al+i\be},&x>0,\\
b\vert x\vert^{\al+i\be},&x\le 0,
\end{array}\right.
\ea
 and their real and imaginary parts.
Here, $a,b\in\CC,\,\vert a\vert+\vert b\vert>0,\,
p\in(0,\iy],\,\be\in\R$, and either $\al>0,\,p\in(0,\iy]$
or $\al=0,\,p=\iy$.
 In particular,
 $F_{\al,\be,1,1}(x)
 =\vert x\vert^{\al+i\be}$ and
 $F_{\al,\be,1,-1}(x)
 =\vert x\vert^{\al+i\be}(\mbox{sgn}\, x)$.
 Note that the proof of Theorem \ref{T1.3} (a) is based
 on the inequality $A_\sa\left(F_{\al,\be,1,-1}\right)_p<\iy$.
 In addition, note that the asymptotic behaviour of
 $E_n\left(F_{\al,\be,a,b}\log^k(\cdot),
 L_p([-1,1])\right),\,k\in\N$, as $n\to\iy$,
 can be found as well by
 using techniques from \cite{B1946.1} and
 \cite[Sect. 8.1]{G2008book}.
\end{remark}

\begin{remark}\label{R1.7}
The proof of Theorem \ref {T1.1}
  follows general ideas developed in \cite[Sect. 2]{G1992},
  \cite[Sect. 11.3]{G2000}, and
  \cite[Corollary 7.6]{G2020},
  and it is based on the contemporary version of the technique
  originally developed in \cite{B1946,G1982}.
  Note that certain elements of this technique are used for
  finding asymptotically sharp constants in Bernstein-Markov-Nikolskii
  type inequalities \cite{GT2017,G2017,G2019b} (see also earlier results
  in \cite{G2005,LL2015a,LL2015b,GM2018}).
\end{remark}

The proofs of Theorems \ref {T1.1} through \ref {T1.4} are
 presented in Section \ref{S3}.
  Section \ref{S2} contains certain properties of functions
and polynomials.

 \section{Properties of Functions  and Polynomials}\label{S2}
 \noindent
\setcounter{equation}{0}
In this section we discuss certain properties of entire functions
 of exponential type, polynomials, and the error $A_\sa(f)_p$
that are needed for the proof of the theorems.
We start with several known properties of univariate  entire functions
 of exponential type.
 \begin{lemma}\label{L2.1}
 The following statements are valid for $\sa>0$:\\
 (a) A crude Nikolskii's inequality
 \bna
   \left\|g\right\|_{L_{\iy}(\R)}
  \le C(p)\sa^{1/p}
  \left\|g\right\|_{L_{p}(\R)},\qquad g\in B_\sa\cap L_p(\R),\quad
 p\in(0,\iy),\label{E2.1}
  \ena
  where $C$ is independent of $g$, holds true.\\
  (b) If $g\in B_\sa\cap L_\iy(\R)$, then
  $\vert g(x+iy)\vert\le \exp[\sa \vert y\vert]\,\|g\|_{L_\iy(\R)},\,
    x\in\R,\,y\in\R$.\\
  (c) For any sequence $\{g_s\}_{s=1}^\iy,\,
g_s\in B_\sa,\,s\in\N,$
with $\sup_{s\in\N}\| g_s\|_{L_\iy(\R)}= C$, there exist a subsequence
$\{g_{s_m}\}_{m=1}^\iy$ and a function $g_0\in B_\sa\cap L_\iy(\R)$
such that equality
\beq\label{E2.1a}
\lim_{m\to\iy}  g_{s_m}=g_0
\eeq
holds true
uniformly on each compact subset of $\R$.\\
(d) Let $\E_\sa$ be the set  of all entire functions
$g(z)=\sum_{k=0}^\iy c_kz^k$,
     satisfying the following condition: for any $\de>0$ there exists
      a constant $C(\de)$,
     independent of $g$ and $k$, such that
     \beq\label{E2.2}
     \vert c_k\vert \le \frac{C(\de)\sa^k(1+\de)^k}{k!},\qquad k=0,\,1,\ldots.
     \eeq
     Then for any sequence $\{g_s\}_{s=1}^\iy\subseteq\E_\sa$ there exist
     a subsequence $\{g_{s_m}\}_{m=1}^\iy$ and a function $g_0\in B_\sa$
     such that equality \eqref{E2.1a}
holds true
     uniformly on each compact subset of $\R$.\\
(e) If $g\in B_0\cap L_\iy(\R)$, then $g$ is a constant function.
  \end{lemma}
  \proof
  The proofs of standard statements (a), (b), and (c)
  can be found, e.g.,  in \cite[Eq. 4.9(29)]{T1963},
  \cite[Lemma 3.4.3]{SW1971}, and \cite[Sect. 3.3.6]{N1969},
  respectively.
  Statement (d) was proved in \cite[Lemma 2.6]{G2017}.

The proof of statement (e),  given in the
  not easily accessible Bernstein's book
  \cite[Corollary 3.8.1]{B1937}, is long and based on
   an estimate of $\vert g^\prime(x)\vert$. For the reader's
    convenience, we present a short proof of (e).
    Indeed, for every $\vep>0,\, g\in B_\vep\cap L_\iy(\R)$.
    Then using statement (b), we obtain
    $\vert g(x+iy)\vert\le \exp[\vep \vert y\vert]\,\|g\|_{L_\iy(\R)},\,
    x\in\R,\,y\in\R$,
    and setting
    $\vep\to 0+$, we see that $g$ is bounded on $\CC$.
    Therefore, statement (e) is valid by Liouville's theorem.
   \hfill $\Box$
  \begin{remark}\label{R2.1a}
  A different proof of Lemma \ref{L2.1} (e) based on
  Bernstein's inequality
    $ \left\|g^\prime\right\|_{L_{\iy}(\R)}
  \le  \vep\left\|g\right\|_{L_{\iy}(\R)},\,
  g\in B_\vep\cap L_\iy(\R)$,
  (see, e.g.,  \cite[Sect. 4.8.2]{T1963}) was suggested by
  one of the referees.
  \end{remark}

  Next, we discuss certain properties of  the error $A_\sa(f)_p$.

  \begin{lemma}\label{L2.2}
  (a) If $A_\sa(f)_p<\iy, \,\sa>0,\,p\in(0,\iy]$, then
  there exists $g\in B_\sa$ such that
  $\|f-g\|_{L_p(\R)}=A_\sa(f)_p$.\\
  (b) For $\sa>0$ and $\g\in\R\setminus\{0\}$,
  \beq\label{E2.4}
  A_\sa(f(\g\cdot))_\iy=A_{\sa/\vert \g\vert}(f)_\iy,
  \eeq
  if either side of \eqref{E2.4} is finite.\\
  (c) If $A_{\sa_0}(f)_p<\iy,\,\sa_0>0,\,p\in(0,\iy]$, then
  $A_{\sa}(f)_p$ is a nonincreasing
   function of $\sa\in [\sa_0,\iy)$
  and
  \beq\label{E2.6}
  A_{\sa-0}(f)_p=A_\sa(f)_p,\qquad\sa> \sa_0,\quad p\in(0,\iy).
  \eeq
  \end{lemma}
  \proof
  (a) This fact is well-known in approximation theory for $f\in L_p(\R)$
  (see, e.g., \cite[Sects. 2.6.2 and 2.6.3]{T1963}). It is easy to extend it to
  any $f$ with $A_\sa(f)_p<\iy$ since there exists $g^*\in B_\sa$ such that
  $f-g^*\in L_p(\R)$ and $A_\sa(f)_p=A_\sa(f-g^*)_p$.

  (b) Assume that $A_{\sa/\vert \g\vert}(f)_\iy<\iy$.
  Then by Lemma \ref{L2.2} (a),
  there exists
  $g\in B_{\sa/\vert \g\vert}$,
  such that
  $
  \left\|f-g\right\|_{L_\iy(\R)}= A_{\sa/\vert \g\vert}(f)_\iy
  $. Next, $g(\g\cdot)\in B_{\sa}$
  and
  \ba
  A_\sa(f(\g\cdot))_\iy\le\left\|f(\g\cdot)-g(\g\cdot)\right\|_{L_\iy(\R)}
  =A_{\sa/\vert \g\vert}(f)_\iy.
  \ea
   The inequality
   $ A_\sa(f(\g\cdot))_\iy\ge A_{\sa/\vert \g\vert}(f)_\iy$
   can be proved similarly.
   Thus \eqref{E2.4} holds true.

  (c) Relation \eqref{E2.6} for $p\in[1,\iy)$
  was proved in \cite[Sect. 99, Lemma 1]{A1965}, and
  for $p\in(0,\iy)$ we prove it similarly.
  Since $A_{\sa_0}(f)_p<\iy$, there exists $g^*\in B_{\sa_0}$ such that
  $f-g^*\in L_p(\R)$ and $A_{\sa_0}(f-g^*)_p=A_{\sa_0}(f)_p$. Therefore,
  without loss of generality we can assume that $f\in L_p(\R)$.

  Given $\sa>\sa_0$,
  there exists
  $g\in B_\sa\cap L_p(R)$ (by Lemma \ref{L2.2} (a))
  such that
  \beq\label{E2.6a}
  \left\|f-g\right\|_{L_p(\R)}= A_\sa(f)_p.
  \eeq
  Next, given $\vep\in (0,\sa/2)$, the function
  $g((1-\vep/\sa)\cdot)\in B_{\sa-\vep}$
  and by
  triangle inequality \eqref{E1.3} and \eqref{E2.6a},
  \bna\label{E2.6b}
  &&A_{\sa-\vep}(f)_p^{\tilde{p}}
  \le \left\|f-g((1-\vep/\sa)\cdot)\right\|_{L_p(\R)}^{\tilde{p}}\nonumber\\
  &&\le \left\|f-g\right\|_{L_p(\R)}^{\tilde{p}}
  +\left\|g(\cdot)-g((1-\vep/\sa)\cdot)\right\|_{L_p(\R)}^{\tilde{p}}
  = A_\sa(f)_p^{\tilde{p}}+I,
  \ena
  where for a given $M>0$,
  \beq\label{E2.6c}
  I
  =\left( \left\|g(\cdot)-g((1-\vep/\sa)\cdot)\right\|_{L_p(-M,M)}^p
  +\left\|g(\cdot)-g((1-\vep/\sa)\cdot)\right\|_
  {L_p(\R\setminus (-M,M))}^p\right)^{\tilde{p}/p}=(I_1+I_2)^{\tilde{p}/p}.
  \eeq
  Since $g\in L_\iy(\R)$ by \eqref{E2.1}, we have by Bernstein's inequality
  (see, e.g., \cite[Sect. 4.8.2]{T1963})
  \beq\label{E2.6d}
  I_1
  \le 2M \left\|g(\cdot)-g((1-\vep/\sa)\cdot)\right\|_{L_\iy(-M,M)}^p
  \le 2M^{p+1}\vep^p\left\|g\right\|_{L_\iy(\R)}^p.
  \eeq
  Further, by triangle inequality \eqref{E1.3},
  \beq\label{E2.6e}
  I_2
  \le \left(2^{\tilde{p}/p}+1\right)^{p/\tilde{p}}
  \left\|g\right\|_{L_p(\R\setminus (-M/2,M/2))}^p.
  \eeq
  Choosing by \eqref{E2.6e} a large enough $M$ such that $I_2<\de^p/2$
   and then choosing by \eqref{E2.6d}
  a small enough $\vep$ such that $I_1<\de^p/2$, we see from
   \eqref{E2.6c} and \eqref{E2.6b} that
   $
    A_{\sa-0}(f)_p^{\tilde{p}}< A_\sa(f)_p^{\tilde{p}}+\de^{\tilde{p}}
    $. Letting $\de\to 0+$, we arrive at \eqref{E2.6}. \hfill $\Box$

  \begin{remark}\label{R2.2a}
  The example of $f(x)=\cos x$ shows that $1=A_{1-0}(f)_\iy>A_{1}(f)_\iy=0$.
  So in general equality \eqref{E2.6} is not valid for $p=\iy$, while
  \ba
  A_{\sa+0}(f)_p=A_\sa(f)_p,\qquad\sa\ge \sa_0,\quad p\in(0,\iy];
  \ea
  for $p\in[1,\iy]$ it was proved in \cite[Sect. 2.6.22]{T1963}
  (see also \cite[Sect. 99, Lemmas 1 and 2]{A1965}),
   and
  for $p\in(0,1)$ it can be proved similarly by using \eqref{E1.3}.
\end{remark}
  One more property of functions from $B_\sa$ is discussed
   in the following lemma.
\begin{lemma}\label{L2.3}
If $f\in M_{K,N}$ and there exists $g\in B_\sa,\,\sa>0$,
 such that $\|f-g\|_{L_p(\R)}=C_3<\iy,\,
p\in(0,\iy]$, then $g\in M_{KC_4,d+N_1}$, where $d:=\lfloor 1/p\rfloor+1,\,
N_1:=\lceil N\rceil$, and $C_4=C_4\left(N,\sa,p,C_3,K\right)$.
\end{lemma}
\proof
The statement is trivial for $p=\iy$ with $C_4=1+C_3/K$, so we assume that
$p\in(0,\iy)$. Let us set
\ba
h(x):=\left(\frac{\sin(x+\theta)}{x+\theta}\right)^{d+N_1},
\ea
where $\theta\in[0,\pi/2]$ and $h\in B_{d+N_1}$.
Then setting $f_1:=fh$, we obtain
\beq\label{E2.7}
\int_\R\left\vert f_1(x)\right\vert^pdx
=\int_\R\left\vert f(x-\theta)\right\vert^p
\left\vert\frac{\sin x}{x}\right\vert^{(d+N_1)p}dx
=\int_{\vert x\vert<1}+\int_{\vert x\vert\ge 1}
=I_1+I_2,
\eeq
where
\bna
&&I_1\le K^p\int_{\vert x\vert<1}(1+\vert x\vert+\theta)^{N_1p}dx
\le 2K^p(2+\pi/2)^{N_1p},\label{E2.8}\\
&&I_2\le 2K^p\int_{1}^\iy\frac{(1+x+\theta)^{N_1p}}{x^{(d+N_1)p}}dx
\le K^p\frac{2(2+\pi/2)^{N_1p}}{dp-1}\label{E2.9}.
\ena
Combining \eqref{E2.7},  \eqref{E2.8}, and \eqref{E2.9}, we have
\beq\label{E2.10}
\left\|f_1\right\|_{L_{p}(\R)}\le KC_5(N,p).
\eeq
Next, setting $g_1:=gh\in B_{d+N_1+\sa}$ and taking account of
the relations
\ba
\left\|f_1-g_1\right\|_{L_{p}(\R)}
\le \left\|f-g\right\|_{L_{p}(\R)}= C_3,
\ea
 we obtain
from triangle inequality \eqref{E1.3} and \eqref{E2.10}
\ba
\left\|g_1\right\|_{L_{p}(\R)}\le \left(C_3^{\tilde{p}}
+\left(C_5K\right)^{\tilde{p}}\right)^{1/\tilde{p}}=KC_6.
\ea
Hence by Nikolskii's inequality \eqref{E2.1},
\beq\label{E2.11}
\left\vert g(x)
\left(\frac{\sin(x+\theta)}{x+\theta}\right)^{d+N_1}\right\vert
\le K CC_6,\qquad \theta\in [0,\pi/2],\quad x\in\R.
\eeq
Then setting $\theta=0$ in \eqref{E2.11}, we see that for
\ba
x\in E_1:=\{x\in\R:\pi/4+k\pi\le x\le 3\pi/4+k\pi,\,
k=0,\,\pm 1,\ldots\}
\ea
the following inequalities hold true:
\beq\label{E2.12}
\left\vert g(x)\right\vert\le KCC_6\left(\frac
{\vert x\vert}{\vert\sin x\vert}\right)^{d+N_1}\le K2^{(d+N_1)/2}CC_6
\vert x\vert^{d+N_1}.
\eeq
Furthermore, setting $\theta=\pi/2$ in \eqref{E2.11}, we see that for
\ba
x\in E_2:=\{x\in\R:-\pi/4+k\pi\le x\le \pi/4+k\pi,\,
k=0,\,\pm 1,\ldots\}\ea
the following relations are valid:
\beq\label{E2.13}
\left\vert g(x)\right\vert
\le KCC_6\left(\frac
{\vert x\vert+\pi/2}{\vert\cos x\vert}\right)^{d+N_1}
\le K2^{(d+N_1)/2}CC_6
(\vert x\vert+\pi/2)^{d+N_1}.
\eeq
Since $E_1\cup E_2=\R$, it follows from \eqref{E2.12} and \eqref{E2.13}
that $g\in M_{KC_4,d+N_1}$, where
\linebreak $C_4=(\pi/\sqrt{2})^{d+N_1}CC_6$.
\hfill $\Box$
\vspace{.12in}\\
For normed rearrangement-invariant spaces a similar result was proved in
\cite[Lemma 3.4]{G1992}. This proof is more elementary than
the one in \cite{G1992}.

In the next two lemmas we discuss estimates of the error of
  polynomial approximation for functions from $B_\sa$.

\begin{lemma}\label{L2.4}
 Let a sequence $\left\{a_n\right\}_{n=1}^\iy$
 of positive numbers
 satisfy the condition
 \beq\label{E2.13a}
 \sup_{n\in\N}\max\{0,a_n- n\}\le C<\iy.
 \eeq
If $g\in B_\sa\cap L_\iy(\R),\,\sa>0$, and $\tau\in(0,1)$,
 then the following inequality holds true:
  \beq\label{E2.14}
  E_n\left(g,L_\iy([-a_n\tau/\sa,a_n\tau/\sa])\right)
  \le C_7(\tau,C)\exp[-C_8(\tau)\,n]\,\|g\|_{L_\iy(\R^1)},
  \qquad n\in\N,
  \eeq
  where
  \beq\label{E2.15}
  C_7(\tau,C)
  :=\frac{2\tau\exp\left[C\sqrt{1-\tau^2}\right]}
  {\sqrt{1-\tau^2}},\qquad
  C_8(\tau):=\log\left(1+\sqrt{1-\tau^2}\right)-\log \tau-\sqrt{1-\tau^2}>0.
  \eeq
  \end{lemma}
  \noindent
  \proof
  It is known (see, e.g., \cite[Sect 5.4.4]{T1963})
  that for any $g\in B_\sa\cap L_\iy(\R),\,n\in\N,\,
  \tau\in(0,1)$, and $\de>0$,
  \ba
   E_n\left(g,L_\iy([-a_n\tau/\sa,a_n\tau/\sa])\right)
   \le \frac{2\exp[a_n\tau \de]}
   {\de\left(\de+\sqrt{1+\de^2}\right)^{n}}
   \|g\|_{L_\iy(\R)}.
   \ea
   Therefore,
   \bna\label{E2.16}
   &&E_n\left(g,L_\iy([-a_n\tau/\la,a_n\tau/\la])\right)\nonumber\\
   &&\le \frac{2\exp\left[C\tau\de\right]}{\de}
   \exp\left[\left(\tau\de
   -\log\left(\de+\sqrt{1+\de^2}\right)\right)\,n\right]
   \|g\|_{L_\iy(\R)}.
   \ena
   Setting $\de=\sqrt{1-\tau^2}/\tau$ in \eqref{E2.16},
   we arrive at \eqref{E2.14} and \eqref{E2.15}.
   \hfill $\Box$

   In case of $a_n=n,\,n\in\N$, versions of  Lemma \ref{L2.4}
  were proved by the author \cite[Lemma 4.1]{G1982}
   and Bernstein \cite[Theorem VI]{B1946}
  (see also \cite[Sect. 5.4.4]{T1963} and
  \cite[Appendix, Sect. 83]{A1965}).  Multivariate inequalities
  like \eqref{E2.14} were obtained
  in \cite{G1982} and \cite{G1991}.

  \begin{lemma}\label{L2.5}
  Let a sequence $\left\{a_n\right\}_{n=1}^\iy$
  of positive numbers
 satisfy condition \eqref{E2.13a} of Lemma \ref{L2.4}.
 Then the following statements
 are valid:\\
 (a) If $g\in  B_\sa\cap M_{K,N},\,\tau\in(0,1)$, and $p\in(0,\iy]$,
  then the following limit relation holds true:
  \beq\label{E2.17}
  \lim_{n\to\iy}E_n\left(g,
  L_p([-a_n\tau/\sa,a_n\tau/\sa])\right)=0.
  \eeq
  (b) If $g\in  B_\sa,\,\tau\in(0,1) $, and $p\in(0,\iy]$,
  and $\mu=1.508879\ldots$ is the positive solution
  to the equation \eqref{E1.8a},
  then the following limit relation  holds true:
  \beq\label{E2.17a}
  \lim_{n\to\iy}E_n\left(g,
  L_p\left(\left[-\frac{a_n\tau}{\mu\sa},\frac{a_n\tau}{\mu\sa}
  \right]\right)\right)=0.
  \eeq
  \end{lemma}
  \proof
  (a) We prove statement (a) in two steps.\vspace{.12in}\\
  \textbf{Step 1.}
  We first assume that $g\in  B_\sa\cap L_2(\R)$.
  By the Paley-Wiener theorem \cite[Theorem 7.2.1]{S1994},
  there exists $\psi\in L_2([-\sa,\sa])$ such that
  \ba
  g(x)
  =(2\pi)^{-1/2}\int_{-\sa}^\sa \psi(t)
  \exp[-it x]\,dt,\qquad x\in\R,
  \ea
  and $\| g\|_{L_2(\R)}=\| \psi\|_{L_2([-\sa,\sa])}$.
Let $P_n\in\PP_n$ be a polynomial of best approximation to
 $\exp[i\cdot]$ in the metric of $L_\iy([-a_n\tau,a_n\tau])$, i.e.,
 \ba
 \left\|\exp[i\cdot]-P_n(\cdot)\right\|_{L_\iy([-a_n\tau/\sa,a_n\tau/\sa])}
 =E_n\left(\exp[i\cdot],L_\iy([-a_n\tau,a_n\tau])\right).
 \ea
   Then the function
 \ba
 Q_n(x):= (2\pi)^{-1/2}\int_{-\sa}^\sa \psi(t)
  P_n(-t x)\,dt,\qquad x\in\R,
  \ea
  is a polynomial from $\PP_n$, and since $\exp[i\cdot]\in B_1$,
  we obtain from Lemma \ref{L2.4}
  \bna\label{E2.18}
  &&E_n\left(g,
  L_\iy([-a_n\tau/\sa,a_n\tau/\sa])\right)
  \le \left\|g-Q_n\right\|_{L_\iy
  ([-a_n\tau/\sa,a_n\tau/\sa])}\nonumber\\
  &&\le (2\pi)^{-1/2} \int_{-\sa}^\sa\vert \psi(t)\vert\,dt
  \sup_{t\in [-\sa,\sa] }\left\|\exp[it\cdot]-P_n(t\cdot)\right\|_{L_\iy
  ([-a_n\tau/\sa,a_n\tau/\sa])}\nonumber\\
  &&\le (\sa/\pi)^{1/2}\|g\|_{L_2(\R)}
  E_n\left(\exp[i\cdot],L_\iy([-a_n\tau,a_n\tau])\right)\nonumber\\
  &&\le (\sa/\pi)^{1/2}C_7(\tau,C)\|g\|_{L_2(\R)}\exp[-C_8(\tau)\,n],
  \ena
  where $C_7$ and $C_8$ are defined by \eqref{E2.15}.
\vspace{.12in}\\
 \textbf{Step 2.}
 Next, let  $g\in  B_\sa\cap M_{K,N}$ and let $N_1:=\lceil N\rceil\ge 0$.
 Then given $\tau\in (0,1)$ and $\vep\in(0,(1-\tau)/(2\tau (N_1+1))]$,
  the function
 \ba
 g_1(z):=g(z)\left[\frac{\sin(\vep \sa z)}{\vep \sa z}\right]^{N_1+1},\qquad z\in \CC,
 \ea
 belongs to $B_{\sa(1+\vep (N_1+1))}\cap L_2(\R)$
 and
 $\left\|g_1\right\|_{L_{2}(\R)}
 \le KC_{9}(\tau,N,\sa)\,\vep^{-N_1-1/2}$.
 Replacing now $g$ with $g_1,\,\sa$ with $\sa_1:=\sa(1+\vep (N_1+1))$, and
 $\tau$ with $\tau_1:=\tau(1+\vep (N_1+1))\le (1+\tau)/2<1$
 in \eqref{E2.18},
 we see from  \eqref{E2.18} that
 \bna\label{E2.19}
 &&E_n\left(g_1,
  L_\iy([-a_n\tau/\sa,a_n\tau/\sa])\right)
  = E_n\left(g_1,
  L_\iy([-a_n\tau_1/\sa_1,a_n\tau_1/\sa_1])\right)\nonumber\\
 &&\le KC_{10}(\tau,N,\sa,C)\vep^{-N_1-1/2}
 \exp\left[-C_8((1+\tau)/2)\,n\right].
 \ena
 Furthermore, using an elementary inequality
 $v-\sin v\le v^3/6,\,v\ge 0$, we have
 \bna\label{E2.20}
 &&\left\|g-g_1\right\|_{L_\iy([-a_n\tau/\sa,a_n\tau/\sa])}\nonumber\\
 &&\le K(N_1+1)(1+a_n/\sa)^N
 \max_{x\in [-a_n\tau/\sa,a_n\tau/\sa]}\left|1-\frac{\sin(\vep\sa x)}
 {\vep\sa x}\right|\nonumber\\
 &&\le (1/6)K(N_1+1)\,\vep^2\,(1+a_n/\sa)^Na_n^2.
 \ena
 Combining \eqref{E1.4}, \eqref{E2.19}, and \eqref{E2.20}, we obtain
 \bna\label{E2.21}
&& E_n\left(g,
  L_\iy([-a_n\tau/\sa,a_n\tau/\sa])\right)\nonumber\\
&& \le K C_{11}(\tau,N,\sa,C))\left(\vep^2\,(n+C)^{N+2}+\vep^{-N_1-1/2}
 \exp[-C_{12}(\tau)\, n]\right).
 \ena
 Then minimizing the right-hand side of \eqref{E2.21} over all
 $\vep\in (0,(1-\tau)/(2\tau (N_1+1))]$,
 we arrive at the following inequality:
\beq\label{E2.22}
E_n\left(g,
  L_\iy([-a_n\tau/\sa,a_n\tau/\sa])\right)
 \le KC_{13}(\tau,N,\sa,C))
 n^{\frac{(N+2)(N_1+1/2)}{N_1+5/2}}
\exp[-C_{14}(\tau,N,C)\, n].
 \eeq
 Finally, \eqref{E2.17} follows from \eqref{E2.22} and the inequality
 \ba
 E_n\left(g,
  L_p([-a_n\tau/\sa,a_n\tau/\sa])\right)
  \le (2(n+C)\tau/\sa)^{1/p}
  E_n\left(g,
  L_\iy([-a_n\tau/\sa,a_n\tau/\sa])\right).
  \ea
  (b) It is known (see, e.g., \cite[Sect 5.4.5]{T1963})
  that for any $g\in B_\sa$ and
  $\tau\in(0,1)$,
  \beq\label{E2.22a}
   E_n\left(g, L_\iy\left(\left[-\frac{a_n\tau}{\mu\sa},\frac{a_n\tau}{\mu\sa}
  \right]\right)\right)
  \le \frac{2\exp[-n\rho]}{1-\exp[-\rho]},
   \eeq
   where
   \ba
   -\rho:=\tau \sqrt{\mu^2+1}/\mu-\log\left( \sqrt{\mu^2+1}+\mu\right)<0,
   \ea
   and  $n\in\N$ satisfies the inequality $n>\tau a_n$.
   The condition on
   $\left\{a_n\right\}_{n=1}^\iy$ guarantees  the existence of
   $n_0(C,\tau)\in\N$
   such that $n>\tau a_n$ for all $n\ge n_0$. Therefore, \eqref{E2.22a}
   holds true for all $n\ge n_0$.
   Then \eqref{E2.17a} follows from \eqref{E2.22a} and the inequality
 \ba
 E_n\left(g, L_p\left(\left[-\frac{a_n\tau}{\mu\sa},\frac{a_n\tau}{\mu\sa}
  \right]\right)\right)
  \le (2(n+C)\tau/(\mu\sa))^{1/p}
  E_n\left(g, L_\iy\left(\left[-\frac{a_n\tau}{\mu\sa},\frac{a_n\tau}{\mu\sa}
  \right]\right)\right).
  \ea
  \hfill $\Box$

  In case of $a_n=n,\,n\in\N$, Lemma \ref{L2.5} (b)
  was proved by   Bernstein \cite{B1947}
  (see also \cite[Sect. 5.4.5]{T1963} and
  \cite[Appendix, Sect. 83]{A1965}).

  We also need a V. A. Markov-type estimate for coefficients of
  a polynomial.
 \begin{lemma}\label{L2.6}
 For any $P(x)=\sum_{k=0}^nc_kx^k,\,n\in\N, \,\vep\in(0,1),\,
 p\in(0,\iy]$, and $k=0,\,1,\,\ldots$,
 the following inequality holds true:
 \beq\label{E2.23}
 \vert c_k\vert \le C_{15}(p,\vep)
  \frac{n^{k+1/p}}{k!(1-\vep)^{k}}\|P\|_{L_p([-1,1])}.
 \eeq
In addition, if a number sequence $\left\{a_n\right\}_{n=1}^\iy$
 satisfies the condition:
 $\sup_{n\in\N}\max\{0,n-a_n\}\le C<\iy$, then there exists
 $n_0=n_0(\vep,C)$
 such that
 \beq\label{E2.24}
 \vert c_k\vert \le C_{16}(p,\vep)
  \frac{\sa^k}{k!(1-\vep)^{2k}}\|P\|_{L_p([-a_n/\sa,a_n/\sa])},\qquad n\ge n_0,\quad
  k=0,\,1,\,\ldots,\quad \sa>0.
 \eeq
\end{lemma}
\proof
Inequality \eqref{E2.23} was proved in \cite[Lemma 2.5]{G2017}.
 Next, applying
\eqref{E2.23} to a polynomial $P(a_n\cdot)$, we see that
\eqref{E2.24} holds true for $C_{16}=C_{15}(1-\vep)^{-1/p}$ and
$n_0=\lfloor C/\vep\rfloor+1$.\hfill $\Box$\\
For $p=\iy$ inequality \eqref{E2.23} is valid for $\vep=0$
 and $C_{15}=1$ (see, e.g., \cite[Eq. 4.8(49)]{T1963}
 with a proof in \cite[Lemma 2.5]{G2019b}).
\vspace{.12in}\\
In the following lemma we discuss a sufficient condition for
$A_\sa(f)_p<\iy$.

\begin{lemma}\label{L2.7}
Let $f\in M_{K,N}$ be a locally integrable function on $\R$
and let us set $d:=\lfloor 1/p\rfloor+1$, where $p\in(0,\iy]$.
In addition, assume that $\F(f)$ is a tempered distribution
and its restriction to $\R\setminus(-a,a),\,a\in (0,\sa)$,
is a $d$-times continuously differentiable function
$h:\R\setminus (-a,a)\to \CC$
such that $h^{(j)}\in L_1(\R\setminus (-a,a)),\,0\le j\le d$.
Then $A_\sa(f)_p<\iy$.
\end{lemma}
\proof
Let $P\in\PP_{2d+1}$ be the Hermite's interpolating polynomial
such that
$P^{(j)}(\pm a)=h^{(j)}(\pm a),\,0\le j\le d$. Then
\ba
h_1(t):=\left\{\begin{array}{ll}
h(t), &t\in \R\setminus (-a,a),\\
P(t), &t\in (-a,a),
\end{array}\right.
\ea
is a $d$-times continuously differentiable function on $\R$ and
$h_1^{(j)}\in L_1(\R),\,0\le j\le d$.
Moreover, integrating by parts, we see that
\ba
\left\vert \F\left(h_1\right)(x)\right\vert
\le (2\pi)^{-1/2}\min\left\{\int_\R\left\vert h_1(t)\right\vert\,dt,
\vert x\vert^{-d}\int_\R\left\vert h_1^{(d)}(t)\right\vert\,dt\right\}
\le C_{17}(f,a,p)(1+\vert x\vert)^{-d},\quad x\in\R.
\ea
Hence $\F\left(h_1\right)\in L_p(\R)$.
Next, the support of the tempered distribution $H:=\F(f)-h_1$
is in $[-a,a]$ since
\ba
\langle\F(f)-h_1,\vphi\rangle
=\int_{\vert t\vert>a}(h(t)-h(t))\vphi(t)\,dt=0
\ea
for all $\vphi\in S(\R)$ with $\mbox{supp}\, \vphi\subseteq
\R\setminus (-a,a)$.
Then by the generalized Paley-Wiener theorem \cite[Theorem 7.2.3]{S1994},
$g=\F(H)$ is an entire function from $B_a$ and
by the Fourier transform inversion formula,
\ba
g(x)=\F\left(\F(f)-h_1\right)(x)=f(-x)-\F(h_1)(x),\qquad x\in\R.
\ea
Hence $f(x)-g(-x)=\F(h_1)(-x)\in L_p(\R)$.
This completes the proof of the lemma.\hfill $\Box$\\
Note that a weaker version of Lemma \ref{L2.7} is proved in
\cite[Lemma 8.3.1]{G2008book}.

   \section{Proofs of Theorems}\label{S3}
 \noindent
\setcounter{equation}{0}
\emph{Proofs of Theorems \ref{T1.1} and \ref{T1.2}.}
We first prove the inequality
\beq\label{E3.1}
\limsup_{n\to\iy}E_n\left(f,L_p([-a_n/\sa,a_n/\sa])\right)
\le A_{\sa-0}(f)_p,\qquad p\in(0,\iy],\quad \sa>\sa_0>0.
\eeq
Assume first that $f\in M_{K,N}$ and $A_{\sa_0}(f)_p<\iy$.
Given $\tau\in(\sa_0/\sa,1)$, there exists
  $g_{\tau\sa}\in B_{\tau\sa}$ (by Lemma \ref{L2.2} (a))
  such that
\beq\label{E3.2}
\left\|f-g_{\tau\sa}\right\|_{L_p(\R)}
=A_{\tau\sa}(f)_p.
\eeq
Then by Lemma \ref{L2.3}, $g_{\tau\sa}\in M_{KC_4,d+N_1}$.
If $f\in L_p(\R),\,p\in(0,\iy]$,  then
 $A_{\sa_0}(f)_p<\iy$
for all $\sa_0>0$.
Hence given $\tau\in(0,1)$, there exists
  $g_{\tau\sa}\in B_{\tau\sa}\cap L_p(\R)$
  (by Lemma \ref{L2.2} (a)),
satisfying \eqref{E3.2}.
By Nikolskii's inequality \eqref{E2.1},
$g_{\tau\sa}\in M_{K,0}$ for some $K>0$.

Therefore, using Lemma \ref{L2.5} (a)
in both cases, we have
\beq\label{E3.3}
\lim_{n\to\iy}
E_n\left(g_{\tau\sa},L_p([-a_n/\sa,a_n/\sa])\right)=0.
\eeq
Next, applying triangle inequality \eqref{E1.4}, we obtain from
\eqref{E3.2} and \eqref{E3.3}
\ba
&&\limsup_{n\to\iy}
E_n^{\tilde{p}}\left(f,L_p([-a_n/\sa,a_n/\sa])\right)\\
&&\le \limsup_{n\to\iy}
E_n^{\tilde{p}}\left(f-g_{\tau\sa},L_p([-a_n/\sa,a_n/\sa])\right)
+\lim_{n\to\iy}E_n^{\tilde{p}}\left(g_{\tau\sa},
L_p([-a_n/\sa,a_n/\sa])\right)\\
&&\le \left\|f-g_{\tau\sa}\right\|_{L_p(\R)}^{\tilde{p}}
= A_{\tau\sa}(f)_p^{\tilde{p}}.
\ea
Then
\beq\label{E3.4}
\limsup_{n\to\iy}
E_n\left(f,L_p([-a_n/\sa,a_n/\sa])\right)
\le A_{\tau\sa}(f)_p.
\eeq
Letting $\tau\to 1-$
in  \eqref{E3.4}, we arrive at \eqref{E3.1}.
Furthermore, note that  the inequalities
\bna
&&\limsup_{n\to\iy}E_n\left(f,L_p([-a_n/\sa,a_n/\sa])\right)
\le A_{\sa}(f)_p,\qquad p\in(0,\iy),\label{E3.5}\\
&&\limsup_{n\to\iy}E_n\left(f,L_\iy([-a_n/\sa,a_n/\sa])\right)
\le A_{\sa-0}(f)_\iy,\label{E3.5a}
\ena
follow from relations \eqref{E3.1} and \eqref{E2.6}.

Thus inequalities \eqref{E3.5} and \eqref{E3.5a} hold true
if either $f\in L_p(\R)$ or
$f\in M_{K,N}$ and $A_{\sa_0}(f)_p<\iy,\,p\in(0,\iy]$.

Next, we prove the
 inequality
 \beq\label{E3.6}
\liminf_{n\to\iy}E_n\left(f,L_p([-a_n/\sa,a_n/\sa])\right)
\ge A_{\sa}(f)_p,\qquad p\in(0,\iy],\quad \sa>\sa_0>0,
\eeq
  by constructing a nontrivial function $g_0\in B_\sa$,
   such that
 \beq \label{E3.7}
 \liminf_{n\to\iy}E_n\left(f,L_p([-a_n/\sa,a_n/\sa])\right)
 \ge \left\|f-g_{0}\right\|_{L_p(\R)}
\ge A_{\sa}(f)_p,\qquad p\in(0,\iy].
 \eeq

 We can assume without loss of generality that $f\in L_p(\R)$.
 Indeed, assume that \eqref{E3.6} is valid for every $f\in L_p(\R)$.
 Then for every $f\in M_{K,N}$ with $A_{\sa_0}(f)_p<\iy$ and
 any $\tau\in(\sa_0/\sa,1)$, there exists $g\in B_{\tau\sa}$
 such that $f-g\in L_p(\R)$. By Lemma \ref{L2.3}, $g\in M_{KC_4,d+N_1}$,
 and using triangle inequality \eqref{E1.4} and
 Lemma \ref{L2.5} (a), we obtain
 \ba
 &&\liminf_{n\to\iy}
E_n^{\tilde{p}}\left(f,L_p([-a_n/\sa,a_n/\sa])\right)\\
&&\ge \limsup_{n\to\iy}
E_n^{\tilde{p}}\left(f-g,L_p([-a_n/\sa,a_n/\sa])\right)
-\lim_{n\to\iy}E_n^{\tilde{p}}\left(g,L_p([-a_n/\sa,a_n/\sa])\right)\\
&&=\limsup_{n\to\iy}
E_n^{\tilde{p}}\left(f-g,L_p([-a_n/\sa,a_n/\sa])\right)
\ge  A^{\tilde{p}}_{\sa}(f-g)_p
=A^{\tilde{p}}_{\sa}(f)_p.
\ea
Thus we can assume that $f\in L_p(\R)$.

 Next, let $\{n_s\}_{s=1}^\iy$ be a subsequence of natural numbers such that
 \beq \label{E3.8}
 \liminf_{n\to\iy}E_n\left(f,L_p([-a_n/\sa,a_n/\sa])\right)
 =\lim_{s\to\iy}E_{n_s}\left(f,L_p
 \left(\left[-a_{n_s}/\sa,a_{n_s}/\sa\right]\right)\right),
 \eeq
 and let $P_n(x)=\sum_{k=0}^nc_{k,n}x^k,\,n\in\N$, be a polynomial,
  satisfying the equality
\beq \label{E3.9}
E_n\left(f,L_p([-a_n/\sa,a_n/\sa])\right)
 =\left\|f-P_n\right\|_{L_p([-a_n/\sa,a_n/\sa])}.
\eeq
It follows from \eqref{E3.9} and \eqref{E1.3} that
\beq \label{E3.10}
\sup_{n\in\N}\left\|P_n\right\|_{L_p([-a_n/\sa,a_n/\sa])}
\le 2^{1/\tilde{p}}\left\|f\right\|_{L_p(\R)}.
\eeq
Then by Lemma \ref{L2.6} and estimate \eqref{E3.10},
\beq\label{E3.11}
 \vert c_{k,n_s}\vert \le 2^{1/\tilde{p}}C_{16}(p,\vep)\left\|
 f\right\|_{L_p(\R)}
  \frac{\sa^k}{k!(1-\vep)^{2k}},\qquad \vep\in(0,1),
  \quad 1\le k\le n_s,
  \quad s\ge n_0^*(\vep),
 \eeq
 for certain $n_0^*\in\N$.
Therefore, \eqref{E3.11} shows that a sequence
of entire functions
$\left\{g_s\right\}_{s=n_0^*}^\iy:=
\left\{P_{n_s}\right\}_{s=n_0^*}^\iy\subseteq\E_\sa$
 satisfies  condition \eqref{E2.2} of Lemma \ref{L2.1} (d) for
 $\de:=(1-\vep)^{-2}-1$.
 Hence there exist
     a subsequence $\{P_{n_{s_m}}\}_{m=1}^\iy$ and a
     function $g_0\in B_\sa$
     such that the relation
     \beq\label{E3.12}
     \lim_{m\to\iy}P_{n_{s_m}}  =g_0
\eeq
     holds true
     uniformly on each compact subset of $\R$.

Then taking account of \eqref{E3.8}, \eqref{E3.9},
and \eqref{E3.12}, we obtain for any $M>0$
\ba
&&\liminf_{n\to\iy}E_n\left(f,L_p([-a_n/\sa,a_n/\sa])\right)
= \lim_{m\to\iy}\left\|f-P_{n_{s_m}}
\right\|_{L_p\left(\left[-a_{n_{s_m}}/\sa,a_{n_{s_m}}/\sa\right]
\right)}\\
&&\ge \lim_{m\to\iy}\left\|f-P_{n_{s_m}}\right\|_{L_p([-M,M])}
=\left\|f-g_0\right\|_{L_p([-M,M])}.
\ea
Hence
\beq\label{E3.13}
\liminf_{n\to\iy}E_n\left(f,L_p([-a_n/\sa,a_n/\sa])\right)
\ge \lim_{M\to\iy}\left\|f-g_0\right\|_{L_p([-M,M])}
=\left\|f-g_0\right\|_{L_p(\R)}
\ge A_\sa(f)_p.
\eeq
Therefore, \eqref{E3.7}  follows from \eqref{E3.13}.

Thus inequality \eqref{E3.6} holds true
if either $f\in L_p(\R)$ or
$f\in M_{K,N}$ and $A_{\sa_0}(f)_p<\iy,\,p\in(0,\iy]$.
Combining \eqref{E3.5}, \eqref{E3.5a}, and \eqref{E3.6},
we arrive at relations
\eqref{E1.17} and \eqref{E1.18} of Theorems
\ref{T1.1} and \ref{T1.2}.\hfill $\Box$\vspace{.12in}\\
\emph{Proof of Theorem \ref{T1.2a}.}
By Lemma \ref{L2.2} (a) there exists $g_0\in B_{\sa_0}$
such that
$A_{\sa_0}(f)_p=\left\|f-g_0\right\|_{L_p(\R)}<\iy,\,p\in(0,\iy]$.
Then by Lemma \ref{L2.5} (b) for $\sa>\mu\sa_0$,
\beq\label{E3.13a}
\lim_{n\to\iy}E_n\left(g_0,
  L_p\left(\left[-{a_n\tau}/{\sa},{a_n\tau}/{\sa}
  \right]\right)\right)
  \le
  \lim_{n\to\iy}E_n\left(g_0,
  L_p\left(\left[-\frac{a_n\tau}{\mu\sa_0},\frac{a_n\tau}{\mu\sa_0}
  \right]\right)\right)=0.
  \eeq
  Since $f-g_0\in L_p(\R)$, we can use Theorem \ref{T1.2}
   for $p\in(0,\iy)$ along with \eqref{E3.13a} and
   triangle inequality \eqref{E1.4}.
   Therefore, we obtain
\ba
&&A_\sa^{\tilde{p}}(f)_p
=A_\sa^{\tilde{p}}\left(f-g_0\right)_p\\
&&=\lim_{n\to\iy}E_n^{\tilde{p}}\left(f-g_0,
  L_p\left(\left[-{a_n\tau}/{\sa},{a_n\tau}/{\sa}
  \right]\right)\right)
  -\lim_{n\to\iy}E_n^{\tilde{p}}\left(g_0,
  L_p\left(\left[-{a_n\tau}/{\sa},{a_n\tau}/{\sa}
  \right]\right)\right)\\
&&\le \lim_{n\to\iy}E_n^{\tilde{p}}\left(f,
  L_p\left(\left[-{a_n\tau}/{\sa},{a_n\tau}/{\sa}
  \right]\right)\right)\\
&&\le \lim_{n\to\iy}E_n^{\tilde{p}}\left(f-g_0,
  L_p\left(\left[-{a_n\tau}/{\sa},{a_n\tau}/{\sa}
  \right]\right)\right)
  +\lim_{n\to\iy}E_n^{\tilde{p}}\left(g_0,
  L_p\left(\left[-{a_n\tau}/{\sa},{a_n\tau}/{\sa}
  \right]\right)\right)\\
&&=A_\sa^{\tilde{p}}\left(f-g_0\right)_p
=A_\sa^{\tilde{p}}(f)_p.
\ea
Thus \eqref{E1.18} for $p\in(0,\iy)$
 is established. Relation \eqref{E1.17}
for $p=\iy$
can be proved similarly.
   \hfill $\Box$\vspace{.12in}\\
\emph{Proof of Theorem \ref{T1.3}.}
If $\al=0,\,2,\,\ldots$ and $\be=0$, then
statements (a) and (b) of Theorem \ref{T1.3}
trivially hold true. So in the proofs of
these statements we assume that
if $\al=0,\,2,\,\ldots$, then $\be\neq 0$.

(a) To prove statement (a), we first show that
$A_\sa(F_{\la,\be,1,-1})_p<\iy,\,\sa>0$,
where $F_{\la,\be,1,-1}=\vert x\vert^{\la+i\be}
(\mbox{sgn}\,x)$
(see also Remark \ref{R1.6}), $\be\in\R,\,p\in(0,\iy],\,\la>0$, and
if $\la=1,\,3,\,\ldots$, then $\be\ne 0$.
Note that $F_{\la,\be,1,-1}$ is continuous on $\R$ and
$F_{\la,\be,1,-1}\in M_{1,\la}$.
It is known \cite[Eqs. (13) and (15), p. 173]{GS1964}
 that under the above condition on $\la$ and $\be$,
 the Fourier transform of the tempered distribution $F_{\la,\be,1,-1}$
 for $x\in \R\setminus \{0\}$ is
\ba
\F\left(F_{\la,\be,1,-1}\right)(x)
=-(2/\pi
)^{1/2}i\cos((\la+i\be) \pi /2)\,\Gamma (\la +i\be+1)\,|x|^{-\la
-i\be-1}(\mbox{sgn}\,x).
\ea
Since $\la>0$, the restriction of $\F\left(F_{\la,\be,1,-1}\right)$
to $\R\setminus(-a,a),\,a\in(0,\sa)$, is
an infinitely differentiable function $h:\R\setminus (-a,a)\to \CC$
such that $h^{(j)}\in L_1(\R\setminus (-a,a)),\,j=0,\,1,\ldots$.
Hence by Lemma \ref{L2.7}, $A_\sa(F_{\la,\be,1,-1})_p<\iy$.
Therefore, for $\la=\al+1$, there exists a function $g\in B_\sa$ such that
$F_{\al+1,\be,1,-1}-g\in L_p(\R),\,\al>-1$. Moreover, setting
$g^*(x):=(g(x)-g(-x))/2\in B_\sa$, we see that $g^*$ is an odd function
and
\ba
\vert x\vert^{\al+1+i\be}
(\mbox{sgn}\,x)-g^*(x)\in L_p(\R),\qquad\al>-1.
\ea
Next, since $f_{\al,\be}\in L_p[-1,1]$ for $\al>-1/p$, it is easy to see that
\ba
\vert x\vert^{\al+i\be}
-{g^*(x)}/{x}\in L_p(\R),\qquad \al>\max\{-1,-1/p\}.
\ea
Thus
\ba
\max\left\{A_\sa\left(f_{\alpha,
\be,c}\right)_p,A_\sa\left(f_{\alpha,\be,s}\right)_p\right\}
\le A_\sa\left(f_{\alpha,\be}\right)_p<\iy,
\quad\al>\max\{-1,-1/p\},\quad p\in(0,\iy],\quad \sa>0,
\ea
 while
statement (a) for $\al=0$ and $p=\iy$ is trivially valid.

Further, we prove the following equalities:
 \beq\label{E3.14}
 A_{\sa-0}\left(f_{\alpha,\be}\right)_\iy
=A_{\sa}\left(f_{\alpha,\be}\right)_\iy,\quad
A_{\sa-0}\left(f_{\alpha,\be,c}\right)_\iy
=A_{\sa}\left(f_{\alpha,\be,c}\right)_\iy,\quad
A_{\sa-0}\left(f_{\alpha,\be,s}\right)_\iy
=A_{\sa}\left(f_{\alpha,\be,s}\right)_\iy,
 \eeq
 that are needed for the proofs of statements (b) and (c).
Using Lemma \ref{L2.2} (b), we see that for $\g>1,\,
A_{\sa/\g}\left(f_{\alpha,\be}\right)_\iy
=\g^\al A_{\sa}\left(f_{\alpha,\be}\right)_\iy$.
Therefore, letting $\g\to 1-$, we see that
the first equality in \eqref{E3.14} holds true.
Next, using Lemma \ref{L2.2} (b) for $\g>1$ again, we obtain
\beq\label{E3.15}
A_{\sa/\g}\left(f_{\alpha,\be,c}\right)_\iy
\le \g^\al\cos(\be\log \g)
A_{\sa}\left(f_{\alpha,\be,c}\right)_\iy
+\g^\al\sin(\be\log \g)
A_{\sa}\left(f_{\alpha,\be,s}\right)_\iy.
\eeq
Letting $\g\to 1-$ in \eqref{E3.15}, we arrive at
the second equality in \eqref{E3.14}.
The third equality in \eqref{E3.14} can be proved similarly.

(b) We first note that by statement (a),
  $A_{\sa_0}\left(f_{\al,\be}\right)<\iy$ for
all $\sa_0>0$ and either
 $\al>\max\{-1,-1/p\},\,
p\in(0,\iy]$ or $\al=0,\,p=\iy$.
In addition, $f_{\al,\be}\in M_{1,\al}$ for
$\al\ge 0$. Therefore, using statements (a) and (b) of
Theorem \ref{T1.1} for $\sa=1$
and $a_n=n,\,n\in\N$, and taking account of
the first equality in \eqref{E3.14},
we obtain
\beq\label{E3.15a}
\lim_{n\to\iy}E_n\left(f_{\al,\be},L_p([-n,n])\right)
= A_{1}(f_{\al,\be})_p,\qquad p\in(0,\iy],\quad \al\ge 0.
\eeq
Then \eqref{E1.19} for $\al\ge 0$ follows immediately
from \eqref{E3.15a} and \eqref{E1.4a} for $\eta=n$.
We cannot use Theorem \ref{T1.1} for $p\in(0,\iy)$ and
$\al\in (\max\{-1,-1/p\},0)$ since
in this case $f_{\al,\be}\notin M_{K,N},\,N\ge 0$.
However, in this case
 we can use Theorem \ref{T1.2a}
 for $\sa=1$
and $a_n=n,\,n\in\N$, since
  by statement (a),
  $A_{\sa_0}\left(f_{\al,\be}\right)_p<\iy$ for
 all $\sa_0>0$. Thus statement (b) is established.

(c) Let us choose a number sequence $\{a_n\}_{n=1}^\iy$  such that
$\lfloor a_n\rfloor=n,\,n\in\N$, which satisfies
the condition of Theorem \ref{T1.1}
with $C=1$. Then using Theorem \ref{T1.1} for
$p\in(0,\iy],\,\al\ge 0,\,\sa=1$, and
Theorem \ref{T1.2a} for
$p\in(0,\iy),\,\al\in (\max\{-1,-1/p\},0),\,\sa=1$,
and taking account of \eqref{E3.14}, we obtain
\bna
\lim_{n\to\iy}E_n\left(f_{\al,\be,c},L_p([-a_n,a_n])\right)
= A_{1}(f_{\al,\be,c})_p,\qquad p\in(0,\iy],\label{E3.15b}\\
\lim_{n\to\iy}E_n\left(f_{\al,\be,s},L_p([-a_n,a_n])\right)
= A_{1}(f_{\al,\be,s})_p,\qquad p\in(0,\iy].\label{E3.15c}
\ena
Next for $\be\neq 0$, let us consider  sequences
\ba
n_k:=\left\lfloor\exp\left[\frac{\pi k}{\be}\right]\right\rfloor,\,
m_k:=\left\lfloor\exp\left[\frac{(2\pi+1) k}{2\be}\right]\right\rfloor,\,
a_{n_k}:=\exp\left[\frac{\pi k}{\be}\right],\,
a_{m_k}:=\exp\left[\frac{(2\pi+1) k}{2\be}\right],\,\, k\in\N.
\ea
Then
\bna
&&a_{n_k}^{-\al}f_{\alpha,\be,c}\left(a_{n_k} x\right)
=a_{m_k}^{-\al}f_{\alpha,\be,s}\left(a_{m_k} x\right)
=(-1)^kf_{\alpha,\be,c}\left(x\right),\qquad k\in\N,
\quad x\in\R,\label{E3.17}\\
&&a_{n_k}^{-\al}f_{\alpha,\be,s}\left(a_{n_k} x\right)
=-a_{m_k}^{-\al}f_{\alpha,\be,c}\left(a_{m_k} x\right)
=(-1)^kf_{\alpha,\be,s}\left(x\right),\qquad k\in\N,
\quad x\in\R.\label{E3.18}
\ena
Finally, relation \eqref{E1.20}
follows from equalities \eqref{E3.15b}, \eqref{E1.4a},
 and \eqref{E3.17},  and \eqref{E1.21}
follows from \eqref{E3.15c}, \eqref{E1.4a}, and \eqref{E3.18}.
\hfill $\Box$\vspace{.12in}\\
\emph{Proof of Theorem \ref{T1.4}.}
Setting $\g=a_{n_k}=\exp[\pi k/\be],\,k\in\N$, in Lemma
\ref{L2.2} (b), we obtain from \eqref{E3.17} and \eqref{E3.18}
that for $\sa>0$,
\bna
&&A_\sa\left(f_{0,\be,c}\right)_\iy
=A_\sa\left(f_{0,\be,c}\left(a_{n_k}\cdot\right)\right)_\iy
=A_{\sa/a_{n_k}}\left(f_{0,\be,c}\right)_\iy,\qquad k\in\N,\label{E3.19}\\
&&A_\sa\left(f_{0,\be,s}\right)_\iy
=A_\sa\left(f_{0,\be,s}\left(a_{n_k}\cdot\right)\right)_\iy
=A_{\sa/a_{n_k}}\left(f_{0,\be,s}\right)_\iy, \qquad k\in\N\label{E3.20}.
\ena
Next, by Lemma \ref{L2.2} (a), there exists
$g_k\in B_{\sa/a_{n_k}}$ such that
$\left\|f_{0,\be,c}-g_k\right\|_{L_\iy(\R)}
=A_{\sa/a_{n_k}}\left(f_{0,\be,c}\right)_\iy$.
Then $\sup_{k\in\N}\|g_k\|_{L_\iy(\R)}\le 2$.
By Lemma \ref{L2.1} (c),
there exist a subsequence
$\{g_{k_m}\}_{m=1}^\iy$ and a function
$g_0\in B_\sa\cap L_\iy(\R)$
such that the equality
$\lim_{m\to\iy}  g_{k_m}=g_0$
holds true
uniformly on each compact subset of $\R$.
Moreover, since  by Lemma \ref{L2.1} (b), for every $\vep>0$ there exists $m_0(\vep)\in\N$
such that for $m\ge m_0$,
$\vert g_{k_m}(x+iy)\vert\le 2\exp[\vep \vert y\vert],\,
    x\in\R,\,y\in\R$, we conclude that
    $g_0\in B_0\cap L_\iy(\R)$. Therefore, $g_0$ is a constant function,
    by Lemma \ref{L2.1} (e).
    Then for any $M>0$,
\ba
\liminf_{m\to\iy}A_{\sa/a_{n_{k_m}}}\left(f_{0,\be,c}\right)_\iy
\ge \lim_{m\to\iy}\left\|f_{0,\be,c}-g_{k_m}\right\|_{L_\iy([-M,M])}
=\left\|f_{0,\be,c}-g_{0}\right\|_{L_\iy([-M,M])}.
\ea
Hence letting $M\to\iy$, we obtain
\beq\label{E3.21}
\liminf_{m\to\iy}A_{\sa/a_{n_{k_m}}}\left(f_{0,\be,c}\right)_\iy
\ge\left\|f_{0,\be,c}-g_{0}\right\|_{L_\iy(\R)}
\ge\inf_{C\in\R}\left\|f_{0,\be,c}-C\right\|_{L_\iy(\R)}.
\eeq
Since
$\inf_{C\in\R}\left\|f_{0,\be,c}-C\right\|_{L_\iy(\R)}
\ge A_\sa\left(f_{0,\be,c}\right)_\iy$,
it follows from \eqref{E3.19} and \eqref{E3.21} that
\ba
A_\sa\left(f_{0,\be,c}\right)_\iy
=\inf_{C\in\R}\left\|f_{0,\be,c}-C\right\|_{L_\iy(\R)}
=\left\|f_{0,\be,c}\right\|_{L_\iy(\R)}
=1.
\ea
Equality
$A_\sa\left(f_{0,\be,s}\right)_\iy=1$ can be proved similarly
if we use \eqref{E3.20} instead of \eqref{E3.19}.
Finally,
\ba
1\ge A_\sa\left(f_{0,\be}\right)_\iy
\ge A_\sa\left(f_{0,\be,c}\right)_\iy=1.
\ea
Thus \eqref{E1.22} is established.\hfill $\Box$
\section*{Acknowledgements}
 We are grateful to the anonymous referees
 for valuable suggestions.

\end{document}